\documentclass{IEEEtran}
\usepackage{amsmath,amssymb,amsthm}
\usepackage{ucs}
\usepackage{hyperref}
\usepackage{cleveref}
\usepackage{autonum}
\usepackage[normalem]{ulem}

\newtheorem{theorem}{Theorem}
\newtheorem{remark}{Remark}
\newtheorem{lemma}[theorem]{Lemma}
\newtheorem{corollary}[theorem]{Corollary}
\usepackage{mathtools}
\DeclarePairedDelimiter{\ceil}{\lceil}{\rceil}
\usepackage{color}
\usepackage[]{algorithm2e}
\usepackage{url}
\usepackage{graphicx}
\usepackage{adjustbox}
\usepackage{graphicx}
\usepackage{multicol}
\usepackage{blindtext}
\usepackage{algpseudocode}
\usepackage{hyperref} 

\begin{document}
\title{Sufficient conditions for the uniqueness of solution of the  weighted norm minimization problem}

\author{K. Z. Najiya, Munnu Sonkar and C. S. Sastry \\
Department of Mathematics \\
Indian Institute of Technology, Hyderabad, 502285, India. \\ Email:\{ma17resch01001, ma17resch11004 and csastry\}@iith.ac.in}
\date{}
\maketitle
\begin{abstract}
Prior support constrained compressed sensing, achieved via the weighted norm minimization, has of late become popular due to its potential for applications. For the  weighted norm minimization problem,   
$$ 
 min \|x\|_{p,w} \text{ subject to } y=Ax, \; p=0,1, \text{ and } w \in [0,1], 
$$
uniqueness results are known when $w=0,1$. Here, $\|x\|_{p,w}=w\|x_T\|_p+\|x_{T^c}\|_p, \; p=0,1$ with $T$ representing the partial support information. The work reported in this paper presents the conditions that ensure the uniqueness of the solution of this problem for general $w \in [0,1]$.  
\end{abstract} 
\section{Introduction}
In Compressed Sensing (CS), a sparse signal $x \in \mathbb{R}^{n}$ can be recovered from a small set of  measurements $y\in \mathbb{R}^{m}$ satisfying $y = Ax$ with $k \ll m$, where $k$ is the number of nonzero elements in $x$. The results that guarantee the uniqueness of the recovery process depend on the restricted isometry property (RIP) of the sensing matrix $A$ \cite{elad2010sparse}\cite{foucart2017mathematical}{\cite{naidu2016deterministic}}. In many applications, one obtains some a priori information about the partial support of the sparse solution to be recovered. 
For instance, in applications involving recovering time-correlated signals \cite{vaswani2010modified}, prior-support constrained sparse recovery attains importance. In  recent years, compressed sensing with a priori support information has caught the attention of several researchers \cite{vaswani2010modified}\cite{friedlander2011recovering}\cite{jacques2010short}\cite{von2007compressed}, to name a few. The weighted norm minimization aims at providing signals, satisfying the data constraint, that are sparse inside and sparsest outside a given prior support. In \cite{vaswani2010modified}, the authors have modified the 1-norm by taking zero weights on the known partial support, minimizing thereby the terms in the complement of prior support set. The results in \cite{vaswani2010modified} have presented the uniqueness of solution of weighted norm minimization under the stated conditions.  When all the weights are set to 1, the weighted 0-norm and the weighted 1-norm problems coincide respectively with their standard 0-norm and 1-norm counterparts, whose exact recovery conditions have been established in \cite{candes2005decoding}.
The authors of \cite{friedlander2011recovering}\cite{liu2014compressed} have established 
the stability of recovery in 
noisy-setting  for weighted 1-norm minimization problem.
To the best of our knowledge, however, the uniqueness of the solution of the general weighted 0-norm and weighted 1-norm minimization problems has not been proposed to date.
Motivated by this, the present work proposes sufficient conditions for the uniqueness of the solution of the weighted 0,1-norm minimization problems. We show that our conditions mostly coincide with those of known cases when the weights are 0, 1.
\par The paper is organized as 6 sections. In sections 2 and 3, we provide basic introduction to Compressed Sensing and existing uniqueness results respectively. In sections 4 and 5, we discuss the uniqueness results with general weights for 0-norm and 1-norm problems respectively. The paper ends with concluding remarks in section 6.  
    


\section{Compressed sensing}
Compressive sensing (CS) \cite{elad2010sparse} is a technique that reconstructs  a signal, which is compressible or sparse in some domain, from a small set of linear measurements.  Let $\sum_{k}^{n} := \{x \in \mathbf{R}^{n}:\|x\|_0 \leq k \}$ be the set of all k-sparse signals in $\mathbf{R}^{n}$. Here $\|x\|_0= |\{ i : x_i \neq 0 \}|$ stands for the number of nonzero components in $x$. For simplicity in notation,  we represent the set $\{1,2, \dots, n\}$ as $[n]$. 
For $A  \in \mathbf{R}^{m \times n}$ with  $m << n$,   suppose  $y=Ax$. One may recover the sparsest solution of this  system from the following minimization problem 
:\begin{equation}\label{P0}   (P_0  ) \; \; min \|x\|_0 \; \; \; \textit{subject to }\;  y= Ax.
\end{equation}
  Since $l_0$ minimization problem becomes  NP-hard as the dimension increases, the convex relaxation of $l_0$ problem has been proposed as
\begin{equation}\label{eq:P1}
   (P_1  ) \; \;  min \|x\|_1 \; \; \;\textit{subject to }\;  y=Ax.
\end{equation}
 The coherence $\mu(A)$ of a matrix $A$ is the largest absolute normalized inner product between different columns of it, that is,
$$
    \mu(A)=  \max_{1\leq i,j \leq n, \; i \neq j}\frac{  |a_i^T a_j|}{\|a_i\|_2 \|a_j\|_2}, 
$$
where $a_i$ denotes the $i$-th column in $A$. 

The $k$-th restricted isometry property ($k$-RIP) constant $\delta_k$ of a matrix $A$ is the smallest real number such that 
$$
(1-\delta _{k})\|x\|_{{2}}^{2}\leq \|Ax\|_{{2}}^{2}\leq (1+\delta _{k})\|x\|_{{2}}^{2}, 
$$
for all $x$ such that $\| x \|_0 \leq k < n$. The restricted orthogonality constant $\theta_{s,\Tilde{s}}$ of a matrix $A$ is the smallest real number such that
$$
    |\eta'A_{T}'A_{\tilde T} \tilde \eta| \leq \theta_{s,\tilde s} \|\eta\|_2 \| \tilde \eta \|_2,  
$$
for all disjoint sets $T$ and $\tilde T$ with $|T| \leq s$ and $| \tilde T |\leq \tilde s $ such that $s+ \tilde{s}\leq n$ and for all vectors $\eta \in \mathbb{R}^{|T|}$ and $\tilde \eta \in \mathbb{R}^{|\tilde T|}$. Here, $A_T$ denotes the restriction of the matrix $A$ to the columns corresponding to the indices in $T \subseteq [n]$. For simplicity, we denote $\theta_{s}:=\theta_{s,s}$ . In \cite{candes2005decoding}, E.  Candes and T. Tao have given the conditions for the exact recovery of $x$ from the pair $(A,y)$ in terms of Restricted Isometry Constant (RIC) for  (\ref{P0}) and (\ref{eq:P1}).
These results, stated in our notation, are as follows:
\begin{theorem}\label{Candes0-thm}(E. Candes et. al.  \cite{candes2005decoding}): Suppose that $s\geq 1$ is such that \begin{equation}\label{condn00}
   \delta_{2s}<1
\end{equation}
and let $N\subseteq [n]$ be such that $|N|\leq s$. Let $y:=Ax$, where $x$ is an arbitrary vector supported on $N$. Then $x$ is the unique minimizer to (\ref{P0}) so that $x$ can be reconstructed from knowledge of the vector $y$ (and $a_i\;'s$).
\qed
\end{theorem}

\begin{theorem}\label{Candes-thm}(E. Candes et. al.  \cite{candes2005decoding}):
Suppose that $s\geq 1$ is such that 
\begin{equation}\label{Cand-condn}
    \delta_s+\theta_{s,s}+\theta_{2s,s} <1
\end{equation}
and let $x$ be a real vector supported on a set $N\subseteq [n]$ obeying $|N|\leq s$. Put $y:=Ax$. Then $x$ is unique minimizer to (\ref{eq:P1}). 
\qed
\end{theorem}
D. Donoho and X. Huo \cite{donoho2001uncertainty} have  shown the exact recovery condition  for $P_1$ in terms of mutual coherence. 
If $x$ is a $k$ sparse vector and matrix $A$ is $k$-RIP compliant, $ k<  \frac{1}{2}\big( 1+\frac{1}{\mu}\big)$ is an exact recovery condition for $P_1$ problem. 
The following result is relevant to the objective of present work.  
\begin{lemma}
\label{Candes-lemma}(E. Candes et. al.  \cite{candes2005decoding}): Let $s\geq 1$ be such that $\delta_s+\theta_{s,2s}<1$, and $c$ be a real vector supported on $N \subseteq [n]$ obeying $|N|\leq s$. Then there exists a vector $\gamma \in  \mathbb{R}^n$ such that $\gamma'a_i=c_i$ for all $i \in N$ where $a_i$ is the $i^{th}$ column of a matrix $A \in \mathbb{R}^ {m\times n}$. Furthermore, $\gamma$ obeys
\begin{equation}\label{Candes-lemmacondn}
|<\gamma,a_i>|\leq \frac{\theta_s}{(1-\delta_s-\theta_{s,2s})\sqrt{s}}.\|c\|, \; \forall i \notin N.
\end{equation}
\qed
\end{lemma}
\section{Compressed sensing with partial support constraint}
It may be noted that the reconstruction method given by $P_1$ in (\ref{eq:P1}) is nonadaptive as no information about $x$ is used in $P_1$. It can, however, be made partially adaptive by imposing constraints on the support of the solution to be obtained.  In \cite{vaswani2010modified}\cite{friedlander2011recovering}\cite{liu2014compressed} (and the references therein) the authors have modified the cost function of $P_1$ problem by incorporating the partial support information into the reconstruction process as detailed below. 
\par Consider that $T$ is the known partial support information of signal $x$. Here $T$ is considered in general sense that it can have an error part which corresponds to the complement of support of $x$. In \cite{vaswani2010modified}, the authors have modified the $P_0$ problem by considering zero weights in $T$ and posed it as follows: 
\begin{equation}  \label{Nama0}
min \|x_{T^c}\|_0 \; \textit{ subject to} \; y=Ax.
\end{equation}
This problem recovers a signal that satisfies the data constraint and whose support is sparsest outside $T$. The following result in \cite{vaswani2010modified} establishes the uniqueness of (\ref{Nama0}).
\begin{theorem}\label{Vaswani0}
(N. Vasawani et. al.  \cite{vaswani2010modified}): 
Given a sparse vector $x$ with support  $N=T\cup\Delta/ \Delta_e$ where $\Delta$ and $T$ are unknown  and known disjoint supports respectively, and $\Delta_e$ is the error in known support such that $ \Delta_e \subseteq T$. Consider reconstructing it from $y=Ax$ by solving (\ref{Nama0}). Then $x$ is the unique minimizer of (\ref{Nama0}) if $\delta_{k+2u}<1$,\noindent where $k:=|T|$ and $u:=|\Delta|$.
\qed
\end{theorem}
\noindent In \cite{vaswani2010modified}, the authors have also considered the convex relaxation of (\ref{Nama0}) as 
\begin{equation}  \label{Nama1}
min \|x_{T^c}\|_1 \; \textit{ subject to} \; y=Ax.
\end{equation}
The uniqueness condition of (\ref{Nama1}) has been established by the following results.

\begin{theorem}(N. Vasawani et. al.  \cite{vaswani2010modified}): \label{NamrataUNIQ}
Given a sparse vector $x$ whose support $N=T \cup \Delta/ \Delta_e$ where $\Delta$ and $T$ are unknown  and 
known disjoint supports respectively, and $\Delta_e$ is the error in known support such that $ \Delta_e \subseteq T$. Consider reconstructing it from $y=Ax$ by solving (\ref{Nama1}). Then $x$ is the unique minimizer of (\ref{Nama1}) if
\begin{enumerate}
    \item $\delta_{k+u}<1$ and $\delta_{2u} +\delta_{k}+ \theta_{k,2u}^2 <1,$
    \item $\rho_k(2u,u)+\rho_k(u,u)<1$, with  $\rho_k(s,\Tilde{s}):=\frac{\theta_{\Tilde{s},s}+\frac{\theta_{\Tilde{s},k}\theta_{s,k}}{1-\delta_k}}{1-\delta_s-\frac{\theta_{s,k}^2}{1-\delta_k}}$,
    \end{enumerate}
\noindent where $s:=|N|$, $k:=|T|$ and $u:=|\Delta|$.
\qed
\end{theorem}
    \begin{corollary} \label{Coro} (N. Vasawani et. al.  \cite{vaswani2010modified}): Given a sparse vector, $x$, whose support $N=T \cup \Delta/ \Delta_e$ where $\Delta$ and $T$ are unknown  and known disjoint supports respectively, and $\Delta_e$ is the error in known support such that $ \Delta_e \subseteq T$. Consider reconstructing it from $y=Ax$ by solving (\ref{Nama1}). Then $x$ is the unique minimizer of (\ref{Nama1}) if $u\leq k$ and $\delta_{k+2u}<\frac{1}{5}$.
    \qed 
    \end{corollary}
\noindent Since sparsity of a signal inside $T$ is unconstrained in (\ref{Nama0}), the recovered signal may not be sparse in $T$. In order to recover a signal, satisfying the data constraint, which is in general sparse inside $T$ and sparsest outside $T$, one may choose general weights $w\in[0,1]$ and propose the general weighted-zero-norm problem:
\begin{equation} \label{weighted0}
    (P_{0,w})\;\;min \|x\|_{0,w} \;\;\; \textit{subject to} \; \;  y=Ax,  
\end{equation}
where $\|x\|_{0,w}=w\|x_T\|_0+\|x_{T^c}\|_0$. It may be noted that when $w=0$, $P_{0,w}$ coincides with (\ref{Nama0}) and when $w=1$, it coincides with the standard $P_0$ problem (\ref{P0}).  As stated in previous section, the uniqueness results in these two cases are established by Theorem \ref{Vaswani0} and Theorem \ref{Candes0-thm} respectively.
In \cite{friedlander2011recovering}, nevertheless, the authors have convexified this problem for  a general weight vector $w \in [0,1]$ and an arbitrary subset $T$ of $[n]$ the following way:
\begin{equation} \label{weighted}
    (P_{1,w})\hspace{.3cm} min \|x\|_{1,w} \; \textit{ subject to } \; y=Ax,
\end{equation}
\noindent where  $\|x\|_{1,w}:=\sum_{i}w_i|x_i|$ with $w_i=  \begin{cases} w \; & \text{for} \; i \in T \\ 1 & \text{ for } i \notin T \end{cases}$  .
\par In general, in applications, $T$ can be drawn  from the estimate of the support of signal or from its largest coefficients.  It has been shown in \cite{friedlander2011recovering} that a signal $x$ can be stably and robustly recovered from $P_{1,w}$ problem in noisy case if  at least $50\%$ of the partial support information is accurate. 
The uniqueness result in Theorem \ref{NamrataUNIQ} holds in a case when $w$ is set to 0 in $P_{1,w}$. In the case, where $w=1$, however, $P_{1,w}$ coincides with $P_{1}$. To the best of our knowledge, the uniqueness of solution of $P_{p,w}$, with $p=0, 1$, is not known for $w \in (0,1)$.  The present work aims at providing the stated uniqueness in the cases complementary to the known cases (viz, $w$ =\,0, 1). 

\section{Uniqueness of solution of weighted 0-norm problem} 
\noindent Our uniqueness result for weighted 0-norm minimization may be summarized in the form of following theorem, which is motivated by the results in \cite{vaswani2010modified}. 
\begin{theorem}\label{thm0} Let $x$ be a real sparse vector supported on $N\subseteq [n]$ with $|N|= s$ and $y=Ax$, where $A\in \mathbb{R}^{m \times n}$ with $m<n$. Let $T \subseteq [n]$, with $|T|= k$ and $\Delta_1=T \cap N$ with $|\Delta_1|= t$ and $\Delta=T^c \cap N$ with $|\Delta|= u$. If 
\begin{equation}\label{condn0}
\delta_{k+2u+\ceil{w t}}<1, 
\end{equation}
then x is the unique minimizer to the $P_{0,w}$ problem in (\ref{weighted0}) for $0\leq w\leq 1$.
\end{theorem}
\noindent \textit{Proof}:
Let $\Tilde{x}$ be a minimizer of (\ref{weighted0}). Then,
$\|\Tilde{x}\|_{0,w}\leq \|x\|_{0,w}$, which implies that 
$\|\Tilde{x}_{T^c}\|_0\leq  w\|x_T\|_0+\|x_{T^c}\|_0- w\|\Tilde{x}_T\|_0\leq w\|x_T\|_0+\|x_{T^c}\|_0\leq w t+u$. Hence,  $\Tilde{x}_{T^c}$ has at most $w t+u$ number of non-zero elements. Therefore $\Tilde{x}$ remains supported on a subset of $T$ of cardinality at most $k$ and and on a set $\Tilde{\Delta}\subseteq T^c$ of cardinality at most $w t+u$. Similarly $x$ is also supported on a subset $\Delta_1\subseteq T$ of cardinality $t\leq k$ and on a set $\Delta\subseteq T^c$ of cardinality at most $u$. Then the support of $\Tilde{x}-x$ remains contained in the union $T\cup\Delta\cup\Tilde{\Delta}$, which is of cardinality at most $k+u+w t+u=k+2u+w t$. Now $A(\Tilde{x}-x)=0$ reduces to $A_{T\cup\Delta\cup\Tilde{\Delta}}(\Tilde{x}-x)=0$. As $0< \delta_{k+2u+\ceil{w t}}<1$,  $A_{T\cup\Delta\cup\Tilde{\Delta}}$ is a full rank matrix, which implies that $ \Tilde{x}=x$.
\qed
\begin{remark}
Here the ceiling operation $\ceil{ w t \;}$ is used to take the smallest integer greater than or equal to the real number $wt$. 
\end{remark}
\begin{remark}
    \par When $w=1$,  the weighted $0$-norm problem coincides with the standard $0$-norm problem in (\ref{P0}) and $k+2u+w t =k-t+2(t+u)= 2s+e$ with $e=|T\cap N^c|$. Further, if  $T\subseteq N$ then $e=0$. Hence $\delta_{k+2u+\ceil{w t}}<1$ coincides with the uniqueness condition $\delta_{2s}<1$ of the standard $0$-norm problem in (\ref{P0}).
    \par When $w=0$, the weighted $0$-norm problem coincides with the $0$-norm problem in  (\ref{Nama0}) and the uniqueness  condition in (\ref{condn0}) of the weighted 0-norm problem  coincides with $\delta_{k+2u}<1$ of Theorem \ref{Vaswani0}.
    
\end{remark}

\section{Uniqueness of solution of weighted 1-norm problem} \label{Present-Work}
\noindent Our uniqueness result for weighted 1-norm minimization is established with the help of following lemma:
\begin{lemma}\label{unique-coro}
Let $x \in \mathbb{R}^n$ be a real sparse vector supported on $N\subseteq[n]$ with $|N|= s$ and $A\in \mathbb{R}^{m \times n}$ with $m<n$ . Let $c \in \mathbb{R}^n$  be such that 
\begin{equation}
c_i =
  \begin{cases}
    w.sgn(x_i) & \text{for~} i \in T \\
    sgn(x_i) & \text{for~} i \in \Delta \\
    0 & \text{otherwise,}
  \end{cases} 
  \end{equation} 
  where $T\subseteq [n]$ with $|T|= k$, $\Delta=T^c\cap N$ with $|\Delta|= u$ and $w\in[0,1]$. If
  \begin{equation}\label{condn11}
   \bigg(\sqrt{\frac{k w^2+u}{k+u}}\bigg)\theta_{k+u} +\delta_{k+u}+\theta_{{k+u},2(k+u)}<1,
\end{equation}
then there exists a vector $\gamma \in  \mathbb{R}^n$ such that 
\begin{enumerate}
    \item $\gamma'a_i=w.sgn (x_i)$ for $i \in T$
    \item $\gamma'a_i=sgn (x_i)$ for $i \in \Delta$
    \item $|\gamma'a_i|< 1 $ for $i \in (T\cup \Delta)^c$. 
\end{enumerate}
\end{lemma}
\noindent \textit{Proof}:
Since $\delta_{k+u}+\theta_{k+u,2(k+u)}<1$ follows from (\ref{condn11}), Lemma \ref{Candes-lemma} implies that there exists a vector $\gamma \in \mathbb{R}^n$ such that $\gamma'a_i=c_i$ for $i \in T\cup \Delta$, that is, $\gamma'a_i=w.sgn(x_i)$ for $i \in T$ and $\gamma'a_i=sgn(x_i)$ for $i \in \Delta$. Again, from  (\ref{Candes-lemmacondn}) and  (\ref{condn11}), we have 
\begin{equation}
\begin{split}
    |\gamma'a_i| &\leq  \frac{ \theta_{k+u}\,\|c\|}{(1-\delta_{k+u}-\theta_{k+u,2(k+u)})\sqrt{k+u})}\\
    &= \frac{\theta_{k+u}(\sqrt{k w^2+u})}{(1-\delta_{k+u}-\theta_{k+u,2(k+u)})\sqrt{k+u}} < 1. \qed  
    \end{split}
\end{equation}
The following result summarizes the uniqueness of solution of weighted 1-norm minimization problem, whose proof is motivated by the results in \cite{candes2005decoding}. 
\begin{theorem}\label{thm} Let $x$ be a real sparse vector supported on $N\subseteq[n]$ with $|N|= s$ and $y=Ax$, where $A\in \mathbb{R}^{m \times n}$ with $m<n$. Let $T \subseteq [n]$ with $|T|=k$ and  $\Delta=T^c \cap N$ with $|\Delta|=u$. If 
\begin{equation}\label{condn}
   \bigg(\sqrt{\frac{k w^2+u}{k+u}}\bigg)\theta_{k+u} +\delta_{k+u}+\theta_{{k+u},2(k+u)}<1,
\end{equation}
then x is the unique minimizer to the $P_{1,w}$ problem in (\ref{weighted}) for  $0\leq w\leq 1$.
\end{theorem}
\noindent \textit{Proof}: By standard convex arguments, there exists one minimizer $\Tilde{x}$ to the problem (\ref{weighted}), which implies that $\|\Tilde{x}\|_{1,w} \leq \|x\|_{1,w}$. Note that $x_i=0$ for $i \in (T\cup N)^c\subseteq N^c$. We have  
 \begin{equation}  \label{p1w-proof}
 \begin{split}
    \|\Tilde{x}\|_{1,w}&=\sum_{i\in T} w|\Tilde{x}_i| +\sum_{i\in T^c}|\Tilde{x}_i|\\
    &=\sum_{i\in T} w|\Tilde{x}_i|+\sum_{i\in \Delta} |\Tilde{x}_i|+\sum_{i\in (T\cup \Delta)^c}|\Tilde{x}_i|\\
    &= \sum_{i\in T} w|x_i+\Tilde{x}_i-x_i|+\sum_{i\in \Delta} |x_i+\Tilde{x}_i-x_i|\\&+\sum_{i\in (T\cup \Delta)^c}|\Tilde{x}_i-x_i|\\
     &\geq \sum_{i\in T} w. sgn(x_i)(x_i+\Tilde{x}_i-x_i)\\
     &+\sum_{i\in \Delta} sgn(x_i)(x_i+\Tilde{x}_i-x_i)+\sum_{i\in (T\cup \Delta)^c}(\Tilde{x}_i-x_i)\\
     &=\sum_{i\in T} w|x_i|+\sum_{i\in \Delta}|x_i|+\sum_{i\in T} w. sgn(x_i)(\Tilde{x}_i-x_i)\\&+\sum_{i\in \Delta} sgn(x_i)(\Tilde{x}_i-x_i)+\sum_{i\in (T\cup \Delta)^c}(\Tilde{x}_i-x_i)\\
     &\geq \|x\|_{1,w}+\sum_{i\in T}\gamma'a_i(\Tilde{x}_i-x_i)+\sum_{i\in \Delta}\gamma'a_i(\Tilde{x}_i-x_i)\\&+\sum_{i\in (T\cup \Delta)^c}\gamma'a_i(\Tilde{x}_i-x_i)\\
     &=\|x\|_{1,w}+\gamma'A(\Tilde{x}-x)
     =\|x\|_{1,w}.
\end{split}
\end{equation}
\noindent In the above chain of steps, the vector $\gamma \in \mathbb{R}^n$ is supposed to satisfy the following properties: 
\begin{enumerate}
    \item $\gamma'a_i=w.sgn (x_i)$ for $i \in T$
    \item $\gamma'a_i=sgn (x_i)$ for $i \in \Delta$
    \item $|\gamma'a_i|< 1 $ for $i \in (T\cup \Delta)^c$. 
\end{enumerate}
In view of (\ref{condn}),  the existence of such a vector $\gamma$ is guaranteed by Lemma \ref{unique-coro}. 
From (\ref{p1w-proof}), it follows that $ \|\tilde x \|_{1,w}=\|x\|_{1,w}$. Consequently, all the inequalities in (\ref{p1w-proof}) must be equalities. But then $\sum_{i\in (T \cup \Delta)^c}|\tilde{x}_i|=\sum_{i \in (T\cup \Delta)^c} (\gamma'a_i)\Tilde{x}_i$ implies that  $\Tilde{x}_i=0$ on $(T\cup \Delta)^c$ as $|\gamma'a_i|< 1$ on $(T\cup \Delta)^c$. Now $Ax=A\Tilde{x}$ reduces to $A_{T\cup\Delta}(x-\Tilde{x})=0$. By (\ref{condn}), we have $\delta_{k+u} <1$  which implies that $\Tilde{x}_i=x_i$ on $T\cup\Delta$. Thus $\Tilde{x}=x$ as claimed.
\qed
\begin{remark}
    When $w=1$,  the weighted $1$-norm problem coincides with the standard $1$-norm problem in  (\ref{eq:P1}) and $k+u=t+u+k-t=s+e$, where $t=|T\cap N|$ and $e=|T\cap N^c|$. Further, if $T\subseteq N$ then $e=0$. In this case, $k+u$  coincides with $s$ and the uniqueness condition  (\ref{condn}) of Theorem \ref{thm} coincides with the uniqueness condition
   $\theta_{s} +\delta_{s}+\theta_{s,2s}<1$
     of the standard $1$-norm problem.
     
    \par When $w=0$, the weighted $1$-norm problem coincides with $1$-norm problem (\ref{Nama1}), and the uniqueness condition gets reduces to\\
    \begin{equation}\label{w0}
         \big(\sqrt{\frac{u}{k+u}}\big)\theta_{k+u} +\delta_{k+u}+\theta_{{k+u},2(k+u)}<1. 
    \end{equation}
   As such, it is not possible to compare the above condition to the uniqueness condition of Theorem \ref{NamrataUNIQ}. This is because, the proofs of both adopt different strategies. 
 In order to deduce a condition from (\ref{w0}) in terms of RIC (that is akin to the condition in Corollary \ref{Coro}), we use the inequality $\theta_{s,\tilde s}\leq \delta_{s+\tilde s}$. Then,  $\theta_{k+u}\leq \delta_{2(k+u)}$ and $\theta_{{k+u},2(k+u)}\leq \delta_{3(k+u)}$. Again if $u\leq k$, then $\frac{u}{k+u}\leq \frac{1}{2}$. Hence,  (\ref{w0}) holds if $(\frac{1}{\sqrt{2}}+2)\delta_{3(k+u)}<1$, that is, $\delta_{3(k+u)}<\frac{\sqrt{2}}{1+2\sqrt{2}}\approx 0.369$.
 
\end{remark}
\section{Conclusion}
The current work has proposed the conditions that guarantee the uniqueness of solution of weighted 0-norm and weighted 1-norm minimization problems for $w \in [0,1]$.  It has been analyzed further that the uniqueness conditions match with their known counterparts in the particular cases where (i). $w=0,1$ with 0-norm, (ii). $w=1$ with 1-norm. In the case where $w=0$ with 1-norm, however, our RIC-condition does not exactly match with its corresponding known condition.  \\

\noindent{\bf Acknowledgments:} The first author is thankful to UGC, Govt. of India, (JRF/2016/409284), for its financial support. The second author is thankful for the support that he receives from MHRD, Government of India. 

\bibliographystyle{plain}

\bibliographystyle{plainnat}
\bibliography{bibliography.bib}
\clearpage

\end{document}